\documentclass[reqno]{amsart}

\usepackage{amsthm, amssymb, latexsym, tabls, hhline,lscape}
\usepackage[enableskew]{youngtab}
\usepackage{diagrams}

\newtheorem{theorem}{Theorem}[section]
\newtheorem{proposition}[theorem]{Proposition}
\newtheorem{corollary}[theorem]{Corollary}

\theoremstyle{definition}
\newtheorem{definition}[theorem]{Definition}
\newtheorem{example}[theorem]{Example}

\theoremstyle{remark}
\newtheorem{remark}[theorem]{Remark}
\numberwithin{equation}{section}
\newcommand{\Sym}{\ensuremath{\mathit{Sym}}}
\newcommand{\QSym}{\ensuremath{\mathit{QSym}}}
\newcommand{\NSym}{\ensuremath{\mathit{NSym}}}
\newcommand{\ra}{\rightarrow}
\newcommand{\ot}{\ensuremath{\otimes}}

\newcommand{\scal}[1]{\langle #1 \rangle}
\newcommand{\shuf}{{\,\sqcup\!\sqcup\,}}

\makeindex

\title[Grothendieck Groups of a Tower of Algebras]
{\bf Algebraic Structures on Grothendieck Groups of a Tower of Algebras}

\author{N. Bergeron}\address[Nantel Bergeron]
{Department of Mathematics and Statistics\\ York  University\\ To\-ron\-to, Ontario M3J 1P3\\ CANADA}
\email{bergeron@mathstat.yorku.ca}
\urladdr{http://www.math.yorku.ca/bergeron}

 \author{H. Li}\address[Huilan Li]
 {Department of Mathematics and Statistics\\ York  University\\ To\-ron\-to, Ontario M3J 1P3\\ CANADA}
 \email{lihuilan@mathstat.yorku.ca}

\date{\today}

\thanks{This work is supported in part by CRC and NSERC.}



\begin{document}
\maketitle

\begin{abstract}
The Grothendieck group of the tower of symmetric group algebras
has a self-dual graded Hopf algebra structure. Inspired by this,
we introduce by way of axioms, a general notion of a tower of
algebras and study two Grothendieck groups on this tower linked by
a natural paring. Using representation theory, we show that our
axioms give a structure of graded Hopf algebras on each
Grothendieck groups and these structures are  dual to each other.
We give some examples to indicate why  these axioms are necessary.
We also give auxiliary results that are helpful to verify the
axioms. We conclude with some remarks on generalized towers of
algebras leading to a structure of generalized bialgebras (in the
sense of Loday) on their Grothendieck groups.
\end{abstract}

\section{Introduction}
\setcounter{equation}{0}

In 1977, L. Geissinger realized \Sym\ (symmetric functions in
countably many variables) as a self-dual graded Hopf
algebra~\cite{LG}. Using the work of Frobenius and
Schur~\cite{AVZ}, this can be interpreted as the self-dual
Grothendieck Hopf algebra of the tower of symmetric group algebras
$\bigoplus_{n\geq0}\mathbb{C}\mathfrak{S}_{n}$. Since then, we
have encountered many instances of combinatorial Hopf algebras. In
each instance, we study a pair of dual Hopf algebras, and find
that this duality can be interpreted as the duality of the
Grothendieck groups of an appropriate tower of algebras. For
example, C. Malvenuto and C. Reutenauer established the duality
between the Hopf algebra of \NSym\ (noncommutative symmetric
functions) and the Hopf algebra of \QSym\ (quasi-symmetric
functions)~\cite{MR}. Later, D. Krob and J. -Y. Thibon showed that
this duality can be interpreted as the duality of the Grothendieck
groups associated with $\bigoplus_{n\geq0}H_{n}(0)$ the tower of
Hecke algebras at q = 0~\cite{KT}. More recently, it was shown
that if one uses $\bigoplus_{n\geq0}HCl_{n}(0)$ the tower of
Hecke-Clifford algebras at q = 0, then one gets a similar
interpretation for the duality between the $Peak$ algebra and its
dual~\cite{BHT}. In~\cite{ANS} Sergeev constructed semi-simple
super algebras Se$_n$ ($n\geq0$) and a characteristic map from the
super modules of Se$=\bigoplus_n$Se$_n$ to Schur's $Q$-functions
$\Gamma=\mathbb{C}[p_1,p_3,\ldots]\subseteq$\Sym. The space
$\Gamma$ is a self-dual graded Hopf subalgebra of \Sym.
In~\cite{HNT} the tower of 0-Ariki-Koike-Shoji algebras
$\bigoplus_{n\geq0}H_{n,r}(0)$ is shown to be related to the
Mantaci-Reutenauer descent algebras \cite{MR2}, and their duals, a
generalization of quasi-symmetric functions, are introduced by
Poirier \cite{Po}.

Our present goal is to describe a general setting that includes
all the examples above. We study the relationship between some
graded algebras $A$ and the algebraic structure on their
Grothendieck groups $G_0(A)$ and $K_0(A)$. More precisely,
$A=(\bigoplus_{n\geq 0}A_{n},\rho)$ is a graded algebra where each
homogeneous component $A_n$ is itself an algebra (with a different
product). We will call $A$ a tower of algebras if it satisfies
some axioms given in Section~\ref{defTower}. This list of axioms
 implies that their Grothendieck groups are graded Hopf
algebras. Moreover, our axioms allow us to define a paring and to show that the
corresponding Grothendieck groups are graded dual to each other. We also discuss
how to weaken our axioms and still get similar results. This is core of
our paper and is  found in Section~3.

In Section 5 we discuss how our axioms may be adapted to verify that the
Grothendieck groups $G_0(A)$ and $K_0(A)$ have a structure of generalized
bialgebra in the sense of Loday~\cite{L}. This leads to the notion of generalized towers of algebras.
 In Section 2 we recall some definitions and
propositions about bialgebras and Grothendieck groups.  In Section 4 we give
some examples. We also give some general results that are helpful to check the axioms.

\section{Notations and Propositions}
We give s brief review of  the theory of bialgebras~\cite{LG}
and Grothendieck groups~\cite{CR} which will be useful for later
discussion.

\begin{definition} Let $K$ be a commutative ring. A $K$-\textit{algebra} $B$ is a
$K$-module with \textit{multiplication}
$\pi: B\otimes_{K}B\rightarrow B$ and \textit{unit map} $\mu:K
\rightarrow B$ satisfying the \textit{associativity} and the
\textit{unitary property}. We denote this algebra by the triple $(B,\pi,\mu)$.

A $K$-\textit{coalgebra} $C$ is a $K$-module with
\textit{comultiplication} $\Delta: C\rightarrow C\otimes C$ and
\textit{counit map} $\epsilon:C\rightarrow R$ satisfying
\textit{coassociativity} and \textit{counitary property}. We denote this coalgebra by the triple
$(C,\Delta,\epsilon)$.

If a $K$-module $B$ is simultaneously an algebra and a coalgebra,
 it is called a \textit{bialgebra} provided these structures are
compatible in the sense that the comultiplication and counit are
algebra homomorphisms. We denote this bialgebra by the 5-tuple
$(B,\pi,\mu,\Delta,\epsilon)$.

A $K$-linear map $\gamma:H\rightarrow H$ on a bialgebra $H$ is an
\textit{antipode} if for all $h$ in $H$, $\Sigma
h_{i}\gamma(h_{i}')=\epsilon(h)1_{H}=\Sigma \gamma(h_{i})h_{i}'$
when $\Delta h=\Sigma h_{i}\otimes h_{i}'$. A \textit{Hopf
algebra} is a bialgebra with an antipode.
\end{definition}

\begin{definition}
An algebra $B$ is a \textit{graded algebra} if there is a direct
sum decomposition $B=\bigoplus_{i\geq 0} B_{i}$ such that
 $\pi(B_{p}\otimes B_{q})\subseteq B_{p+q}$,
and $\mu(K)\subseteq B_{0}$.

A coalgebra $C$ is a \textit{graded coalgebra} if there is a
direct sum decomposition $C=\bigoplus_{i\geq 0} C_{i}$ such that
$\Delta(C_{n})\subseteq\bigoplus(C_{k}\otimes C_{n-k})$ and
$\epsilon(C_{n})=0$ if $n\geq 1$.

A bialgebra $H=\bigoplus_{n\geq0} H_{n}$ over $K$ is called
\textit{graded connected} if  $H_0=K1_H$. It is well known that a
graded connected bialgebra is a Hopf algebra~\cite{DM}.
\end{definition}

For  $H=\bigoplus_{n\geq0} H_{n}$ a graded bialgebra, its
\textit{graded dual} $H^{*gr}=\bigoplus_{n\geq0} H_n^*$ is also a
graded bialgebra if all $H_n$ are finitely generated and
$H_i^*\otimes H_k^*\cong(H_i\otimes H_k)^*$ for all $i$ and $k$.

We now recall the  definition of Grothendieck groups.
Let $B$ be an arbitrary algebra. Denote
        $$\begin{array}{rcl}
         _{B}\mbox{mod}&=&\mbox{the category of all finitely generated left  $B$-modules,}\\
          \mathcal{P}(B)&=&\mbox{the category of all finitely generated projective left $B$-modules.}
       \end{array}$$

\begin{definition} Let $\mathcal{C}$ be one of the above categories. Let
$\mathbf{F}$ be the free abelian group generated by the symbol
$(M)$, one for each isomorphism class of modules $M$ in
$\mathcal{C}$. Let $\mathbf{F_{0}}$ be the subgroup of
$\mathbf{F}$ generated by all expressions $(M)-(L)-(N)$
arising from all short exact sequences
              $$0\rightarrow L\rightarrow M\rightarrow N\rightarrow0$$
in $\mathcal{C}$. The
\textit{Grothendieck group} ${\mathcal K}_{0}(\mathcal{C})$ of the category
$\mathcal{C}$ is defined by the quotient $\mathbf{F}/\mathbf{F_{0}}$.
For $M\in\mathcal{C}$, we denote by $[M]$ its image in ${\mathcal K}_{0}(\mathcal{C})$.
We then set
             $$G_{0}(B)={\mathcal K}_{0}(_{B}\mbox{mod}) \quad\hbox{and}\quad
            K_{0}(B)={\mathcal K}_{0}(\mathcal{P}(B)).$$
\end{definition}

Now let $B$ be a finite-dimensional algebra over a field $K$. Let
$\{V_{1},\cdots, V_{s}\}$ be a complete list of nonisomorphic
simple $B$-modules. Then their projective covers $\{P_{1},\cdots,
P_{s}\}$ is a complete list of nonisomorphic indecomposable
projective $B$-modules~\cite{ARS}. With these lists, we have

\begin{proposition}
          $$G_{0}(B)=\bigoplus^s_{i=1}\mathbb{Z}[V_{i}]
              \quad\quad\hbox{and}\quad\quad
            K_{0}(B)=\bigoplus^s_{i=1}\mathbb{Z}[P_{i}].$$
\end{proposition}

Let $A$ be an algebra and $B\subseteq A$ a subalgebra, or more generally
let $\varphi\colon B\to A$ be an injection of algebra preserving unities. Let $M$ be
a (left) $A$-module and $N$ a (left) $B$-module. Then the
\textit{induction} of $N$ from $B$ to $A$ is
Ind$^{A}_{B}N=A\otimes_{B}N$ a (left) $A$-module and the
\textit{restriction} of $M$ from $A$ to $B$ is
Res$^{A}_{B}M=\mbox{Hom}_{A}(A,M)$ a (left) $B$-module. In the case of
$\varphi\colon B\to A$, the expression  $A\otimes_{B}N$ is the tensor
$A\otimes N$ modulo the relations $a\otimes bn\equiv a\varphi(b)\otimes n$,
and the left $B$-action on $\mbox{Hom}_{A}(A,M)$ is defined by $bf(a)=f(a\varphi(b))$,
for $f\in \mbox{Hom}_{A}(A,M)$ and $b\in B$.

\section{Grothendieck Groups of a Tower of Algebras}
We now present our axiomatic definition of a {\sl tower} of
algebras. The starting ingredient is a graded algebra
$A=(\bigoplus_{n\geq 0}A_{n},\rho),$ such that each homogeneous
component is itself a finite dimensional algebra. For all $n,m\ge
0$, we require the maps $\rho_{n,m}$ obtained from the products
$\rho$ restricted to $A_n\otimes A_m$ to be injective
homomorphisms of algebras (preserving unities). Our axioms will
allow us to define a notion of induction and restriction on the
Grothendieck groups $G_0(A)=\bigoplus_{n\geq 0}G_0(A_{n})$ and
$K_0(A)=\bigoplus_{n\geq 0}K_0(A_{n})$. This will be the basic
construction to put a structure of graded dual Hopf algebras on
$G_0(A)$ and $K_0(A)$.

\subsection{Tower of Algebras (Preserving unities)} \label{defTower}
Let $A=(\bigoplus_{n\geq 0}A_{n},\rho)$ be a graded algebra. We call it \textit{a tower of
algebras} over field $K=\mathbb{C}$ if the following conditions
are satisfied:
\smallskip

\noindent(1) $A_{n}$ is a finite-dimensional algebra with unit
$1_n$, for each $n$. $A_{0}\cong K$.
\medskip

\noindent(2) The (external) multiplication $\rho_{m,n}:
A_{m}\otimes A_{n}\rightarrow A_{m+n}$ is an injective
homomorphism of algebras, for all $m\mbox{ and }n$ (sending
$1_{m}\otimes 1_{n}$ to $1_{m+n}$ ).
\medskip

\noindent(3) $A_{m+n}$ is a two-sided projective $A_{m}\otimes
A_{n}$-module with the action defined by $a\cdot(b\ot
c)=a\rho_{m,n}(b\ot c)\mbox{ and }(b\ot c)\cdot a=\rho_{m,n}(b\ot
c)a$, for all $m,n\geq0,\ a\in A_{m+n},\ b\in A_{m},\ c\in
A_{n}\mbox{ and }m,n\geq0$.
\medskip

\noindent(4) For every primitive idempotent $g$ in $A_{m+n}$,
$A_{m+n}g\cong\bigoplus (A_m\ot A_n)(e\ot f)$ as (left) $A_m\ot
A_n$-modules if and only if $gA_{m+n}\cong\bigoplus (e\ot
f)(A_m\ot A_n)$ as (right) $A_m\ot A_n$-modules for the same index
of idempotents $(e\ot f)$'s in $A_m\ot A_n$.
\medskip

\noindent(5) The following equality holds for $G_0(A)=\bigoplus_{n\geq 0}G_0(A_{n})$
(or equivalently for $K_0(A)=\bigoplus_{n\geq 0}K_0(A_{n})$)
        $$\begin{array}{ccl}
            &&[\mbox{Res}_{A_k\ot A_{m+n-k}}^{A_{m+n}}\mbox{Ind}_{A_m\ot A_n}^{A_{m+n}}(M\ot N)]\\
           &&\quad = \ \sum_{t+s=k}[\widetilde{\mbox{Ind}}_{A_t\ot A_{m-t}\ot A_s\ot A_{n-s}}^{A_{k}\ot
                A_{m+n-k}}(\mbox{Res}_{A_t\ot A_{m-t}}^{A_m}M\ot\mbox{Res}_{A_s\ot A_{n-s}}^{A_n}N)]
           \end{array}$$
for all $0<k<m+n$, $M$ an $A_{m}$-module and
 $N$ an $A_{n}$-module, or $M$ a projective $A_{m}$-module and
 $N$ a projective $A_{n}$-module. Here
the twisted induction
       $$\begin{array}{ccl}
       &&\widetilde{\mbox{Ind}}_{A_t\ot A_{m-t}\ot A_s\ot
        A_{n-s}}^{A_{k}\ot A_{m+n-k}}(M_1\ot M_2)\ot(N_1\ot N_2)\\
        &&\quad = \ (A_{k}\ot A_{m+n-k})\widetilde{\bigotimes}_{A_{t}\ot A_{m-t}\ot
          A_{s}\ot A_{n-s}}((M_1\ot M_2)\ot(N_1\ot N_2)).\\
         \end{array}$$
This is the usual tensor quotient by the (twisted) relations
        $$\begin{array}{ccl}
        &&(a\ot b)\ot[(c_1\ot c_{2})\cdot(w_1\ot w_2)\ot(d_1\ot d_{2})\cdot(u_1\ot u_2)]\\
        &&\qquad\equiv\ [a\rho_{t,s}(c_1\ot d_1)\ot b\rho_{m-t,n-s}(c_{2}\ot d_{2})]\ot(w_1\ot
        u_1\ot w_2\ot u_2).
       \end{array}$$

Condition (1) guarantees that their Grothendieck groups
 are graded connected. Conditions (2) and (3) insure that  the induction
 and restriction are well defined on $G_0(A)$ and $K_0(A)$.
The duality follows from (4). Finally (5) gives an analogue of
Mackey's formula. This gives us the compatibility relation between
the multiplication and comultiplication that we will define on
$G_0(A)$ and $K_0(A)$.

\subsection{Induction and Restriction on $G_0(A)$}
For  $M$ a left $A_{m}$-module and $N$ a left $A_{n}$-module,
let $M\ot N$ be the left
$A_{m}\ot A_{n}$-module defined by $(a\ot b)\cdot(w\ot
u )=aw\ot bu$, for $a\in A_{m},\, b\in A_{n},\,w\in M \mbox{ and
}u\in N. $ We define
induction  on $G_0(A)$ as follows:

$$\begin{array}{ccccl}
i_{m,n}& : &G_{0}(A_{m})\bigotimes G_{0}(A_{n})
& \rightarrow & G_{0}(A_{m+n}) \\
& & [M]\otimes [N]& \mapsto & [\mbox{Ind}^{A_{m+n}}_{A_{m}\otimes
A_{n}}M\otimes N].
\end{array}$$
For the restriction, we define
$$\begin{array}{ccccc}
 r_{k,l}& : & G_{0}(A_{k+l})& \rightarrow & G_{0}(A_{k})\bigotimes
G_{0}(A_{l})\\
 & & [N]& \mapsto & [\mbox{Res}^{A_{k+l}}_{A_{k}\otimes A_{l}}N].
\end{array}$$

\begin{proposition}\label{ir}$i$ and $r$ are well defined on
$G_{0}(A)$.\end{proposition}

\begin{proof} If $d:M_1\ra M_2$ and $\delta:N_1\ra N_2$ are
isomorphisms, then
    $$A_{m+n}\bigotimes_{A_m\ot A_n}(M_1\ot N_1)\cong A_{m+n}\bigotimes_{A_m\ot A_n}(M_2\ot N_2)
     $$
with the map
  $$(1\ot_{A_{m+n}} (d\ot\delta))(a\ot w\ot u)\stackrel{\rm def}{=}a\ot (d(w)\ot\delta(u)).$$
This is well defined since
$$ (1\ot_{A_{m+n}} (d\ot\delta))(a\ot (bw\ot cu)) =  a\ot (d(bw)\ot\delta(cu))
     =  a\ot (bd(w)\ot c\delta(u)) $$
  $$\qquad =  a\rho(b\ot c)\ot (d(w)\ot \delta(u))
     =  (1\ot_{A_{m+n}} (d\ot\delta))(a\rho(b\ot c)\ot (w\ot u)).
$$
Hence $[\mbox{Ind}^{A_{m+n}}_{A_{m}\otimes A_{n}}M_1\otimes
N_1]=[\mbox{Ind}^{A_{m+n}}_{A_{m}\otimes A_{n}}M_2\otimes N_2]$.
Without loss of generality, assume $[M]=[M']+[M'']$. So there is
a short exact sequence
        $$0\ra M'\ra M\ra M''\ra0.$$
Since $N$ is a finitely generated left $A_{n}$-module, it is a
projective $K$-module. We have
        $$0\rightarrow M'\ot  N\rightarrow M\ot  N\rightarrow M''\ot
N\rightarrow0$$ exact as $K$-modules (also exact as
$A_{m}\otimes A_{n}$-modules). Since $A_{m+n}$ is a (right)
projective $A_{m}\otimes A_{n}$-module, we have
\begin{multline*}0\ra
A_{m+n}\ot_{A_{m}\ot A_{n}} (M'\ot N)\ra A_{m+n}\ot_{A_{m}\ot
A_{n}}
(M\ot N)\\
\ra A_{m+n} \ot_{A_{m}\ot A_{n}} (M''\ot
N)\rightarrow0\end{multline*}exact. Hence
        $$[\mbox{Ind}^{A_{m+n}}_{A_{m}\otimes A_{n}}M\otimes
N]=[\mbox{Ind}^{A_{m+n}}_{A_{m}\otimes A_{n}}M'\otimes
N]+[\mbox{Ind}^{A_{m+n}}_{A_{m}\otimes A_{n}}M''\otimes N].$$
Similarly,
        $$[\mbox{Ind}^{A_{m+n}}_{A_{m}\otimes A_{n}}M\otimes
N]=[\mbox{Ind}^{A_{m+n}}_{A_{m}\otimes A_{n}}M\otimes
N']+[\mbox{Ind}^{A_{m+n}}_{A_{m}\otimes A_{n}}M\otimes N'']$$ for
$[N]=[ N']+[N'']$. Hence $i$ is well defined on $G_{0}(A)$
by induction.

Now we show that $r$ is well defined. Given that
$\mbox{Hom}_{A_{n}}(A_{n},M)\cong M$ for any $A_{n}$-module
$M$ we have that if $N_1\cong N_2$ then $\mbox{Hom}_{A_{n}}(A_{n},N_1)\cong
N_1\cong N_2\cong\mbox{Hom}_{A_{n}}(A_{n},N_2)$. That is
$[\mbox{Res}^{A_{n}}_{A_{k}\ot A_{l}}N_1]=
[\mbox{Res}^{A_{n}}_{A_{k}\ot A_{l}}N_2]$. Without loss of
generality, assume $[N]=[N']+[N'']$. So there is a short exact
sequence
        $$0\ra N'\ra N\ra N''\ra0$$
of $A_n$-modules. Then we have
        $$0\ra \mbox{Hom}_{A_{n}}(A_{n},N')\ra \mbox{Hom}_{A_{n}}(A_{n},N)
          \ra \mbox{Hom}_{A_{n}}(A_{n},N'')\ra0$$
exact as $A_{k}\ot A_{l}$-modules. Hence
        $$[\mbox{Res}^{A_{n}}_{A_{k}\ot A_{l}}N]= [\mbox{Res}^{A_{n}}_{A_{k}\ot
        A_{l}}N']+[\mbox{Res}^{A_{n}}_{A_{k}\ot A_{l}}N'']$$
and again $r$ is well defined by induction on $G_0(A). $\end{proof}

We can now define a multiplication and a comultiplication using $i$
and $r$ and define a unit and a counit on $G_{0}(A)$ as follows:
$$\begin{array}{ccl}
\pi & : & G_{0}(A)\bigotimes
G_{0}(A)\rightarrow G_{0}(A) \\
&& \mbox{where}\
\pi|_{G_{0}(A_{k})\bigotimes G_{0}(A_{l})}=i_{k,l}\\
\Delta & : & G_{0}(A)\rightarrow G_{0}(A)\bigotimes
G_{0}(A) \\
&&
\mbox{where}\ \Delta|_{G_{0}(A_{n})}=\sum_{k+l=n}r_{k,l}\\
\mu & : & \mathbb{Z}\ra G_{0}(A)\\
&& \mbox{where }\mu(a)=a[K]\in G_{0}(A_{0}),\mbox{ for }a\in \mathbb{Z}\\
\epsilon&:& G_{0}(A)\ra \mathbb{Z}\\
&& \mbox{where }\epsilon([M])=\left\{\begin{array}{ll}
                       a & \mbox{if }[M]=a[K],\mbox{ where }a\in \mathbb{Z},\\
                       0 & \mbox{otherwise.}
                       \end{array}
                \right.
\end{array}$$

In Section~3.5 we will prove the associativity of $\pi$, the unity
of $\mu$, the coassociativity of $\Delta$ and the
counity of $\epsilon$.

\subsection{Induction and Restriction on $K_0(A)$}

As before, we define induction and restriction on $K_0(A)$:

$$\begin{array}{ccccl}
i'_{m,n}& : &K_{0}(A_{m})\bigotimes  K_{0}(A_{n})&
\rightarrow & K_{0}(A_{m+n}) \\
& & [P]\otimes [Q]& \mapsto & [\mbox{Ind}^{A_{m+n}}_{A_{m}\otimes
A_{n}}P\otimes Q],
\end{array}$$
and
$$\begin{array}{ccccc}
 r'_{k,l}& : & K_{0}(A_{k+l})& \rightarrow & K_{0}(A_{k})\bigotimes
K_{0}(A_{l})\\
 & & [R]& \mapsto & [\mbox{Res}^{A_{k+l}}_{A_{k}\otimes A_{l}}R].
\end{array}$$

\begin{proposition}$i'$ and $r'$ are well defined on
$K_{0}(A)$.\end{proposition}

\begin{proof}  The proof is similar to Proposition \ref{ir} we
only need to show here that $\mbox{Ind}^{A_{m+n}}_{A_{m}\otimes
A_{n}}P\otimes Q=A_{m+n}\bigotimes_{A_{m}\otimes A_{n}}(P\otimes
Q)$ is a projective $A_{m+n}$-module. Assume that $P\oplus
P'\cong(A_{m})^s \ \mbox{and }Q\oplus Q'\cong (A_{n})^t$ for some
$s$ and $t$. Since
     $$\hbox{$ A_{m+n}\otimes_{A_{m}\otimes A_{n}}(P\otimes Q)\bigoplus
         A_{m+n}\otimes_{A_{m}\otimes A_{n}}(P'\otimes Q)
         \bigoplus
         A_{m+n}\otimes_{A_{m}\otimes A_{n}}((A_{m})^s\otimes Q')$}$$
      $$\begin{array}{cll}
      \quad  &\cong \quad A_{m+n}\bigotimes_{A_{m}\otimes A_{n}}((A_{m})^s\otimes (A_{n})^t)
                  &\cong \quad A_{m+n}\bigotimes_{A_{m}\otimes A_{n}}(A_{m}\otimes A_{n})^{st}\\
      \quad  &\cong \quad (A_{m+n}\bigotimes_{A_{m}\otimes A_{n}}(A_{m}\otimes A_{n}))^{st}
                  & \cong \quad (A_{m+n})^{st},
       \end{array}$$
we have that  $\mbox{Ind}^{A_{m+n}}_{A_{m}\otimes A_{n}}P\otimes Q$ is a
projective $A_{m+n}$-module.

Assume $R\oplus R'\cong (A_{n})^s$ for some $s$. Then there is
a split short exact sequence
        $$0\ra R\ra (A_{n})^s \ra R'\ra 0.$$
Since $\mbox{Hom}_{A_{n}}(A_{n},M)\cong M$ as $A_k\ot
A_l$-modules for any $k+l=n$ and any $A_{n}$-module $M$, the
short sequence
        $$0\ra \mbox{Hom}_{A_{n}}(A_{n},R)\ra \mbox{Hom}_{A_{n}}(A_{n},(A_{n})^s)
\ra \mbox{Hom}_{A_{n}}(A_{n},R')\ra 0$$ is exact and split. That
means
        $$\mbox{Hom}_{A_{n}}(A_{n},(A_{n})^s) \cong \mbox{Hom}_{A_{n}}(A_{n},R)\oplus
         \mbox{Hom}_{A_{n}}(A_{n},R')$$
as $A_{k}\ot A_{l}$-modules. Since
$\mbox{Hom}_{A_{n}}(A_{n},A_{n}^s)\cong A_{n}^s$ and $A_{n}$
is a projective (left) $A_{k}\ot A_{l}$-module, so is
$(A_{n})^s$. {}From above, it follows that
        $$\mbox{Hom}_{A_{n}}(A_{n},R)\oplus
\mbox{Hom}_{A_{n}}(A_{n},R')\cong(A_{n})^s$$ as $A_{k}\ot
A_{l}$-modules, i.e., $\mbox{Hom}_{A_{n}}(A_{n},R)$ is a
summand of $(A_{n})^s$. Therefore,\break
$\mbox{Hom}_{A_{n}}(A_{n},R)$ is a projective $A_{k}\ot
A_{l}$-module. \end{proof}

Using $i'$ and $r'$ we also define a multiplication, a
comultiplication, a unit and a counit on $K_0(A)$.
$$\begin{array}{ccl}
\pi' & : & K_{0}(A)\bigotimes
K_{0}(A)\rightarrow K_{0}(A) \\
& & \mbox{where}\
\pi'|_{K_{0}(A_{k})\bigotimes K_{0}(A_{l})}=i'_{k,l}\\
\Delta' & : & K_{0}(A)\rightarrow K_{0}(A)\bigotimes
K_{0}(A) \\
& &
\mbox{where}\ \Delta'|_{K_{0}(A_{n})}=\sum_{k+l=n}r'_{k,l}\\
\mu' & : & \mathbb{Z}\ra K_{0}(A)\\
&& \mbox{where }\mu'(a)=a[K]\in K_{0}(A_{0}),\mbox{ for }a\in \mathbb{Z}\\
\epsilon' &:& K_{0}(A)\ra \mathbb{Z}\\
&& \mbox{where }\epsilon'([M])=\left\{\begin{array}{ll}
                       a & \mbox{if }[M]=a[K],\mbox{ where }a\in \mathbb{Z},\\
                       0 & \mbox{otherwise.}
                       \end{array}
                \right.
\end{array}$$

In Section 3.5, we will see that the operations above have the desired properties.

\subsection{Pairing on $K_{0}(A)\times G_{0}(A)$}

To show the duality between $K_{0}(A)$ and
$G_{0}(A)$ we define a pairing
 $\scal{\,\,,\,}: K_{0}(A)\times G_{0}(A)\ra \mathbb{\mathbb{Z}} $ where
        $$\scal{[P],[M]}=\left\{\begin{array}{ll}
                       \mbox{dim}_{K}\big(\mbox{Hom}_{A_{n}}(P,M)\big)
                       & \mbox{if }[P]\in K_{0}(A_{n})
                       \mbox{ and }[M]\in G_{0}(A_{n}),\\
                       0 & \mbox{otherwise.}
                       \end{array}
                \right.$$
We also define $\scal{\,\,,\,}: (K_{0}(A)\ot K_{0}(A))\times
(G_{0}(A)\ot G_{0}(A))\ra \mathbb{\mathbb{Z}} $ where
        $$\begin{array}{l}
            \scal{[P]\ot[Q],[M]\ot[N]}\\
            \quad=\left\{\begin{array}{ll}   \mbox{dim}_{K}\big(\mbox{Hom}_{A_{k}\ot A_{l}}
                       (P\ot Q,M\ot N)\big) &
                       \mbox{if }[P]\ot[Q]\in K_{0}(A_{k})\ot
                       K_{0}(A_{l})\\
                         &\mbox{ and }[M]\ot [N]\in G_{0}(A_{k})\ot G_{0}(A_{l}),\\
                       0 & \mbox{otherwise.}
                       \end{array}
                \right.
          \end{array}$$

\begin{proposition}\label{eqs} $\scal{\,\,,\,}$ is a well-defined bilinear pairing on $K_{0}(A)\times
G_{0}(A)$ satisfying the following identities:
        $$\begin{array}{rcl}
\scal{[P]\ot[Q],[M]\ot[N]}&=&\scal{[P],[M]}\scal{[Q],[N]},\\
\scal{\pi'([P]\ot [Q]),[M]}&=&\scal{[P]\ot[Q],\Delta[M]},\\
\scal{\Delta' [P],[M]\ot[N]}&=&\scal{[P],\pi([M]\ot[N])},\\
\scal{\mu'(1),[M]}&=&\epsilon([M]),\\
\scal{[P],\mu(1)}&=&\epsilon'([P]).
\end{array}$$\end{proposition}

\begin{proof}
Assume $[M]=[M']+[M'']$. We have $0\ra M'\ra M\ra M''\ra 0$ a short exact sequence.
Since $P$ is a projective module, the short sequence
     $$ 0\ra\mbox{Hom}_{A_{n}}(P,M')\ra\mbox{Hom}_{A_{n}}(P,M) \ra\mbox{Hom}_{A_{n}}(P,M'')\ra 0$$
is exact. Hence $\scal{[P],[M]}=\scal{[P],[M']}+\scal{[P],[M'']}$.
If we assume  $[P]=[P']+[P'']$, then  $P\cong P'\oplus P''$ and we have
$\mbox{Hom}_{A_{n}}(P,M) \cong \mbox{Hom}_{A_{n}}(P',M)\oplus\mbox{Hom}_{A_{n}}(P'',M)$.
Hence  $\scal{[P],[M]}= \scal{[P'],[M]}+\scal{[P''],[M]}.$ Therefore $\scal{\,\,,\,}$
is a well-defined bilinear pairing on $K_{0}(A)\times G_{0}(A)$.

The identity $\scal{[P]\ot[Q],[M]\ot[N]}=\scal{[P],[M]}\scal{[Q],[N]}$ is clear from the isomorphism
$ \mbox{Hom}_{A_{k}\ot A_{l}}(P\ot  Q,M\ot  N)\cong \mbox{Hom}_{A_{k}}(P,M)\otimes_{K}
 \mbox{Hom}_{A_{l}}(Q,N)$.
For $\scal{\pi'([P]\ot [Q]),[M]}=\scal{[P]\ot[Q],\Delta[M]}$ we use the  Adjointness Theorem~\cite{CR}. We have
        $$\mbox{Hom}_{A_{k+l}}(\mbox{Ind}^{A_{k+l}}_{A_{k}\ot A_{l}}P\ot Q,M)
            \cong\mbox{Hom}_{A_{k+l}}(A_{k+l}\ot_{A_{k}\ot A_{l}}(P\ot Q),M)$$
        $$\begin{array}{rl}
            &\cong\mbox{Hom}_{A_{k}\ot A_{l}}(P\ot Q,\mbox{Hom}_{A_{k+l}}(A_{k+l},M))\\
            &\cong\mbox{Hom}_{A_{k}\ot A_{l}}(P\ot Q,\mbox{Res}^{A_{k+l}}_{A_{k}\ot A_{l}}M),
            \end{array}$$
which gives us
       $$\mbox{dim}_{K}\big(\mbox{Hom}_{A_{k+l}}(\mbox{Ind}^{A_{k+l}}_{A_{k}\ot
           A_{l}}P\ot Q,M)\big)=\mbox{dim}_{K}\big(\mbox{Hom}_{A_{k}\ot
            A_{l}}(P\ot Q,\mbox{Res}^{A_{k+l}}_{A_{k}\ot A_{l}}M)\big).$$
Thus $\scal{i'_{k,l}([P]\ot [Q]),[M]}=\scal{[P]\ot[Q],r_{k,l}[M]}$ and the desired identity follows.

To show $\scal{\Delta'[P],[M]\ot[N]}=\scal{[P],\pi([M]\ot[N])}$, we
need to prove the identity
$\scal{r'_{k,l}[P],[M]\ot[N]}=\scal{[P],i_{k,l}([M]\ot[N])}$, for all
$[P]\in K_{0}(A_{k+l}),\ [M]\in K_{0}(A_{k})\mbox{ and } [N]\in
G_{0}(A_{l})$. This is not as straightforward as before. Here we need the equality
\begin{eqnarray}
    && \mbox{dim}_K\big(\mbox{Hom}_{A_{k+l}}(P,A_{k+l}\ot_{A_{k}
     \ot A_{l}}(M\ot N))\big)\nonumber\\
     \label{dim} &&\qquad\quad = \mbox{dim}_K\big(\mbox{Hom}_{A_{k}\ot
A_{l}}(\mbox{Hom}_{A_{k+l}}(A_{k+l},P),M\ot N)\big).
\end{eqnarray}
Clearly, without lost of generality  we can restrict our attention
to indecomposable projective modules $P$. For such a $P$, there is
a primitive idempotent $g\in A_{k+l}$ such that $P\cong A_{k+l}g$.
We know that for any finite-dimensional algebra $B$ over $K$, $M$
a left $B$-module and $e$ a primitive idempotent, we have
$\mbox{\rm Hom}_B(Be,M)\cong eM$ as vector spaces (see~\cite{NT}).
Hence
      $$\mbox{dim}_{K}\big(\mbox{Hom}_{A_{k+l}}(A_{k+l}g,A_{k+l}\ot_{A_{k}\ot
         A_{l}}(M\ot N))\big) = \mbox{dim}_{K}\big(gA_{k+l}\ot_{A_{k}\ot A_{l}}(M\ot N)\big).
      $$
To prove (\ref{dim}), we expect that
       $$\mbox{dim}_{K} \big(gA_{k+l}\ot_{A_{k}\ot A_{l}}(M\ot N)\big)=\mbox{dim}_{K}\big
           (\mbox{Hom}_{A_{k}\ot A_{l}}(A_{k+l}g\downarrow_{A_k\ot A_l},M\ot N)\big).$$
Since $gA_{k+l}\cong\bigoplus (e\ot f)(A_k\ot A_l)$ as a (right) $A_k\ot
A_l$-module for some idempotents $(e\ot f)$'s in $A_k\ot A_l$, we have
        $$\begin{array}{rl} gA_{k+l}\ot_{A_{k}\ot A_{l}}(M\ot N)
            &\cong\bigoplus (e\ot f)(A_k\ot A_l)\ot_{A_{k}\ot A_{l}}(M\ot N)\\
           & \cong\bigoplus (e\ot f)(M\ot N).\end{array}$$
At the same time from condition (4) $A_{k+l}g\cong\bigoplus
(A_k\ot A_l)(e\ot f)$ as a (left) $A_k\ot A_l$-module for the same
idempotents $(e\ot f)$'s in $A_k\ot A_l$. Hence
        $$\begin{array}{rl}
              \mbox{Hom}_{A_{k}\ot A_{l}}(A_{k+l}g\downarrow_{A_k\ot A_l},M\ot N)
            &\cong\quad \mbox{Hom}_{A_{k}\ot A_{l}}(\bigoplus (A_k\ot A_l)(e\ot f),M\ot N)\\
            &\cong\quad\bigoplus (e\ot f)(M\ot N).\\
             \end{array}$$
 Therefore (\ref{dim}) holds.

We know $\mu'(1)=[K]$ and
        $$\mbox{dim}_{K}\big(\mbox{Hom}_{K}(K,M)\big)=\left\{\begin{array}{ll}
                       a & \mbox{if }[M]=a[K],\mbox{ where }a\in \mathbb{Z},\\
                       0 & \mbox{otherwise,}
                       \end{array}
                \right.$$ therefore
$\scal{\mu'(1),[M]}=\epsilon([M])$. Similarly,
$\scal{[P],\mu(1)}=\epsilon([P])$.\end{proof}

Let $\{V_{1},\cdots, V_{s}\}$ be a complete list of nonisomorphic simple $A_{n}$-modules. Then the
set of their projective covers $\{P_{1},\cdots, P_{s}\}$ is a
complete list of nonisomorphic indecomposable projective
$A_{n}$-modules.
The proposition below is well known (see \cite{CR}).

\begin{proposition}\label{dual}
  $\scal{[P_{i}],[V_{j}]}=\delta_{i,j}$ \ for $1\leq i,j\leq s$.
\end{proposition}

\subsection{Main Result 1}

\begin{theorem} \label{m1}(1) $\pi$ and $\pi'$ are associative. Hence  $\big(G_{0}(A),\ \pi,\ \mu\big)$
and $\big(K_{0}(A),\ \pi',\ \mu'\big)$ are algebras.

\noindent(2) $\Delta$ and $\Delta'$ are coassociative. Hence $\big(G_{0}(A),\
\Delta,\ \epsilon\big)$ and $\big(K_{0}(A),\ \Delta',\ \epsilon'\big)$ are coalgebra.

\noindent(3) If $G_0(A)$ satisfies the condition (5), then
$\Delta$ and $\epsilon$ are algebra homomorphisms and $G_{0}(A)$
is a connected graded  Hopf algebra, as is $K_{0}(A)$ by duality.
Equivalently, the same results holds if instead $K_0(A)$ satisfies
the  condition (5).
\end{theorem}

\begin{proof} (1) We only need to show the associativity of $\pi$. The associativity of
$\pi'$ follows from Prop.~\ref{eqs} and~\ref{dual}. From the associativity of $\rho$ it
is straightforward to verify that
     $$\begin{array}{rl}
     & \mbox{Ind}_{A_{l+m}\otimes A_{n}}^{A_{l+m+n}}(\mbox{Ind}_{A_{l}\otimes
           A_{m}}^{A_{l+m}}L\otimes M)\otimes N \\
     &\qquad\quad =  A_{l+m+n}\bigotimes_{A_{l+m}\otimes A_{n}}((A_{l+m}\bigotimes_{A_{l}\otimes
        A_{m}}(L\otimes M))\otimes N)\\
     &\qquad\quad =  A_{l+m+n}\bigotimes_{A_{l}\otimes A_{m}\otimes A_{n}}(L\otimes M\otimes N)\\
     &\qquad\quad = A_{l+m+n}\bigotimes_{A_{l}\otimes A_{l+n}}(L\otimes(A_{m+n}
        \bigotimes_{A_{m}\otimes A_{n}}(M\otimes N)))\\
     &\qquad\quad = \mbox{Ind}_{A_{l}\otimes A_{m+n}}^{A_{l+m+n}}L\otimes(\mbox{Ind}_{A_{m}\otimes
       A_{n}}^{A_{m+n}} M\otimes N).
    \end{array} $$
Hence $i_{l+m,n}\cdot (i_{l,m}\otimes 1_{n})=i_{l,m+n}\cdot
(1_{l}\otimes i_{m,n})$ and the associativity of $\pi$ follows.

(2) Again we only need to show the coassociativity of $\Delta$, that is,
$(r_{l,m}\otimes 1 )\cdot r_{l+m,n}= (1\otimes r_{m,n})\cdot
r_{l,m+n}$. {}From the definition of $r$ we have
    $$\begin{array}{rl}
      \mbox{Res}_{A_{l}\otimes A_{m}\otimes A_{n} }^{A_{l+m}\otimes A_{n}}\mbox{Res}_{A_{l+m}\otimes
      A_{n}}^{A_{l+m+n}}V
     &=\mbox{Hom}_{A_{l+m}\ot A_{n}}(A_{l+m}\ot A_{n},\mbox{Hom}_{A_{l+m+n}}(A_{l+m+n},V))\\
     &\cong \mbox{Hom}_{A_{l+m+n}}(A_{l+m+n}\ot_{A_{l+m}\ot A_{n}}(A_{l+m}\ot A_{n}),V)
     \end{array}$$
and
        $$\begin{array}{rl}
         \mbox{Res}_{A_{l}\otimes A_{m}\otimes A_{n} }^{A_{l}\otimes A_{m+n}}\mbox{Res}_{A_{l}\otimes
         A_{m+n}}^{A_{l+m+n}} V
        &=\mbox{Hom}_{A_{l}\ot A_{m+n}}(A_{l}\ot A_{m+n},\mbox{Hom}_{A_{l+m+n}}(A_{l+m+n},V))\\
        &\cong \mbox{Hom}_{A_{l+m+n}}(A_{l+m+n}\ot_{A_{l}\ot A_{m+n}}( A_{l}\ot A_{m+n}),V).
       \end{array}$$
Now we want to show that
\begin{eqnarray}
    & \mbox{Hom}_{A_{l+m+n}}(A_{l+m+n}\ot_{A_{l+m}\ot
       A_{n}}(A_{l+m}\ot A_{n}),V)\nonumber\\
   \label{h1}
    &\qquad\quad \cong \mbox{Hom}_{A_{l+m+n}}(A_{l+m+n}\ot_{A_{l}\ot
     A_{m+n}}( A_{l}\ot A_{m+n}),V)
\end{eqnarray}
as $A_{l}\ot A_{m}\ot A_{n}$-modules. In fact,
$A_{l+m+n}\ot_{A_{l+m}\ot A_{n}}(A_{l+m}\ot A_{n})\cong A_{l+m+n} $
as $ A_{l+m+n}\mbox{-}(A_{l+m}\ot A_{n})\mbox{-bimodules}$ and
 $A_{l+m+n}\ot_{A_{l}\ot A_{m+n}}( A_{l}\ot A_{m+n})\cong A_{l+m+n}$
as $A_{l+m+n}\mbox{-}(A_{l}\ot  A_{m+n})\mbox{-bimodules}$. Hence (\ref{h1}) holds
as $A_{l}\ot A_{m}\ot A_{n}$-modules with the action
defined by
$((a\ot b\ot c)\cdot f)(d)=f(d\rho_{l,m,n}(a\ot b\ot c))$
for $a\in A_{l}$, $b\in A_{m}$, $c\in A_{n}$, $d\in A_{l+m+n}$ and
$f\in \mbox{Hom}_{A_{l+m+n}}(A_{l+m+n},V)$. This completes the
proof.

(3) Without loss of generality, we suppose $G_0(A)$ satisfies the
identity in condition (5). For $[M]\in G_{0}(A_{m}),\ [N]\in
G_{0}(A_{n})$, we know that
        $$\Delta(\pi([M]\ot[N]))=\sum^{m+n}_{k=0}\big[\mbox{Hom}_{A_{m+n}}
(A_{m+n},A_{m+n}\ot_{A_{m}\ot A_{n}}(M\ot N))\downarrow_{A_{k}\ot
A_{m+n-k}}\big].$$
We use ``$\downarrow_{A_{k}\ot A_{m+n-k}}$'' to remind us that the module should be
viewed as an $A_{k}\ot A_{m+n-k}$-module. On the other hand, we have in  $A\otimes A$  the following product
        $$\begin{array}{l}
        \Delta[M]\Delta[N]
         = \sum^{m+n}_{k=0}\sum_{t+s=k} \\
        \qquad \big[(A_{k}\ot A_{m+n-k})
          \hbox{$\widetilde{\bigotimes}_{A_{t}\ot A_{m-t}\ot A_{s}\ot A_{n-s}}$}
           \big(\mbox{Hom}_{A_{m}}(A_{m},M)\ot\mbox{Hom}_{A_{n}}(A_{n},N)\big)\big].
           \end{array}$$

To prove that $\Delta$ is an algebra homomorphism we need
$\Delta(\pi([M]\ot[N]))=\Delta[M]\Delta[N]$. For this it is enough to show the equality
of the corresponding terms  for all $0\leq k\leq m+n$ in the expressions above.
When $k=0$, $A_{0}\cong K$ we have
        $$\begin{array}{l}
    \big[\mbox{Hom}_{A_{m+n}}(A_{m+n},A_{m+n}\ot_{A_{m}\ot A_{n}}(M\ot N))\downarrow_{A_{0}
        \ot A_{m+n}}\big]\\
    \qquad= \big[\mbox{Hom}_{A_{m+n}}(A_{m+n},A_{m+n}\ot_{A_{m}\ot A_{n}}(M\ot N))\big]\\
    \qquad= \big[A_{m+n}\ot_{A_{m}\ot A_{n}}(M\ot N)\big]\\
    \qquad= \big[(A_{0}\ot A_{m+n})\widetilde{\bigotimes}_{A_{0}\ot A_{m}\ot A_{0}\ot
        A_{n}}(M\ot N)\big]\\
    \qquad= \big[(A_{0}\ot A_{m+n})\widetilde{\bigotimes}_{A_{0}\ot A_{m}\ot A_{0}\ot
       A_{n}}\big(\mbox{Hom}_{A_{m}}(A_{m},M)\ot\mbox{Hom}_{A_{n}}(A_{n},N)\big)\big].
\end{array}$$
A similar computation holds for $k=m+n$. For $0<k<m+n$, the equality follows from our condition (5):
        $$\begin{array}{l}
         [\mbox{Res}_{A_k\ot A_{m+n-k}}^{A_{m+n}}\mbox{Ind}_{A_m\ot A_n} ^{A_{m+n}}(M\ot N)]\\
      \qquad\quad = \sum_{t+s=k}[\widetilde{\mbox{Ind}}_{A_t\ot A_{m-t}\ot A_s\ot A_{n-s}}^{A_{k}
          \ot A_{m+n-k}}(\mbox{Res}_{A_t\ot A_{m-t}}^{A_m}M\ot\mbox{Res}_{A_s\ot A_{n-s}}^{A_n}N)]
       \end{array}$$
We have that $(G_{0},\pi,\mu,\Delta,\epsilon)$ is a graded bialgebra, hence a graded Hopf algebra.
By duality $K_{0}(A)$ is also a graded Hopf algebra.
\end{proof}

Now we are in the position to state our first main result:

\begin{theorem}If $A$ is a tower of algebras satisfying conditions (1)-(5),
then we can construct on their Grothendieck groups $G_0(A)$ and
$K_0(A)$ a bialgebra structure as above. Moreover,
$(G_{0},\pi,\mu,\Delta,\epsilon)$ and
$(K_{0},\pi',\mu',\Delta',\epsilon')$ are dual to each other as
connected graded bialgebras.\end{theorem}

\subsection{Tower of Algebras (not Preserving unities) and Result 2}
In~\cite{BHRZ}, we consider a semi-tower of algebras with $\rho$
not preserving unities. If we weaken the condition of $\rho$ and
modify the definitions of induction and restriction we can still
get results similar as above. We include only a sketch of the
ideas; the details can be found in~\cite{BL}.

 The usual definitions of induction and restriction as in Section~2 may cause problems
 when $\rho$ does not preserve the unities. For this we need to find a weaker definition.
Let $\varphi\colon B\to A$ be an algebra injection not necessarily
preserving unities. Let $M$ be a left $A$-module. We let
Res$^{A}_{B}M=\{x\in M\ |\ \varphi(1)x=x\}\subseteq M$ be a
submodule. For $x\in\hbox{Res}^{A}_{B}M$ and $b\in B$ the action
is defined by  $\varphi(b)x$. When $\varphi$ preserves the
unities, clearly this definition agrees with the one in Section~2.
For induction, we have to be careful only in the case of
projective modules. Assume that $P$ is an indecomposable left
$B$-modules (we extend our definition linearly). Hence $P\cong Be$
for some primitive idempotent $e\in B$. We let
Ind$^{A}_{B}P=A\varphi(e)$. Again, when $\varphi$ preserves the
unities, it is straightforward to check that this agrees with the
definition of induction in Section~2.

With this in hand, one can adapt all the steps in Section~3.2,~3.3 and~3.4 to obtain

\begin{theorem}[\cite{BL}] If $A$ is a tower of algebras satisfying conditions (1)-(5),
then we can construct on their Grothendieck groups $G_0(A)$ and
$K_0(A)$ a bialgebra structure as above. Moreover,
$(G_{0},\pi,\mu,\Delta,\epsilon)$ and
$(K_{0},\pi',\mu',\Delta',\epsilon')$ are dual to each other as
graded bialgebras.\end{theorem}

\section{Examples}
In this section, we verify that
$\bigoplus_{n\geq0}\mathbb{C}\mathfrak{S}_n$ and
$\bigoplus_{n\geq0}H_{n}(0)$ satisfy all the axioms listed in
Section~\ref{defTower}. They are towers of algebras and we already know that
their Grothendieck groups are dual Hopf algebras, respectively.
We also give an example of graded algebra which do
not satisfy all the axioms  and consequently,
its Grothendieck groups are not dual Hopf algebras.

\subsection{Examples Satisfying All  the Axioms}

\begin{example} Let $A=(\bigoplus_{n\geq0}A_n,\rho)$ with $A_n=\mathbb{C}\mathfrak{S}_{n}$,
where $\mathfrak{S}_n$ is the $n$-permutation group, and
          $$\rho_{m,n}:\mathbb{C}\mathfrak{S}_{m}\otimes\mathbb{C}
            \mathfrak{S}_{n}\rightarrow\mathbb{C}\mathfrak{S}_{m+n},$$
where $\rho_{m,n}(\sigma\otimes\tau)=\sigma(1)\sigma(2)
\cdots\sigma(m)(m+\tau(1))(m+\tau(2))\cdots(m+\tau(n))$. We use
the one line notation of permutations. For example,
$\rho_{2,3}(21\otimes312)=21534.$ It is clear that $\rho$
preserves unities and satisfies associativity. It is also easy to
check that $\rho$ is injective and preserves multiplication.

Now since $\mathbb{C}\mathfrak{S}_{n}$ is a semi-simple algebra, we
know that $\mathbb{C}\mathfrak{S}_{m+n}$ is a two-sided projective
$\mathbb{C}\mathfrak{S}_{m}\otimes\mathbb{C}\mathfrak{S}_{n}$-module.

For finite group $G$,
simple left modules are obtained from primitive idempotents $g\in \mathbb{C}G$.
It is easy to show that the left module $Gg$ is isomorphic to the right module $gG$
(look at their characters). The condition (4) for $A=(\bigoplus_{n\geq0}A_n,\rho)$ is thus satisfied.
Condition (5) is just the Mackey Theorem~\cite{AVZ}.

Hence $A=\bigoplus_{n\geq0}\mathbb{C}\mathfrak{S}_{n}$ is a tower
of algebras and since $\mathbb{C}\mathfrak{S}_{n}$ is a
semi-simple algebra we have that the Grothendieck group
$G_0(A)=K_0(A)$ is a self-dual graded Hopf algebra. The
characteristic map $\hbox{ch}\colon G_0(A)\to\Lambda$, where
$\hbox{ch}([V]=
\sum_{\mu}z_{\mu}^{-1}\hbox{tr}X_{\mu}^{V}p_{\mu}$, is then an
isomorphism of graded Hopf algebras between $G_0(A)$ and $\Lambda$
the Hopf algebra of symmetric functions (see \cite{Ma}).
\end{example}

\begin{remark}
The Sergeev algebra $\hbox{Se}_n$ is the cross product of
symmetric group $\mathfrak{S}_n$ and the Clifford algebra
$\hbox{Cliff}_n$~\cite{ANS}, which is a semisimple superalgebra.
Here consider the Grothendieck groups in categories of finitely
generated supermodules and finitely generated projective
supermodules over these superalgebras. One can modify our
axioms to sit in the category of supermodules over superalgebras.
Its Grothendieck groups $G_0$ and $K_0$ coincide and
have the Hopf algebra structure which is self-dual. It is possible
to check that this tower satisfies the modified
conditions (1)-(5). And $\bigoplus_{n\geq0}\hbox{Se}_n$ is a tower
of superalgebras.
\end{remark}

\begin{example}
Let $A=(\bigoplus_{n\geq0}H_{n}(0),\rho)$ be the direct sum of Hecke
algebras~\cite{KT} where $\rho$ is defined by
$\rho_{m,n}(T_i\otimes 1)=T_i$ and $\rho_{m,n}(1\otimes
T_j)=T_{j+m}$. The $T_i$ for $1\leq i\leq m-1$ are the generators of
$H_m(0)$  satisfying
$$\begin{array}{rclc}
T_i^2&=&-T_i,&\\
T_iT_j&=&T_jT_i& |i-j|>1,\\
T_iT_{i+1}T_i&=&T_{i+1}T_iT_{i+1}.&
\end{array}$$
It is easy to check that $\rho$ preserves unities and satisfies
associativity. Since the $T_i$'s satisfy the braid relations, one
can associate to each permutation $\sigma\in \mathfrak{S}_n$ the
element $T_{\sigma}$ in $H_n(0)$ defined by
$T_{\sigma}=T_{i_1}\cdots T_{i_r}$, where $s_{i_1}\cdots s_{i_r}$
is an arbitrary reduced decomposition of $\sigma$ and $s_i$ is
the simple transposition $(i,\,\,i+1)$. The set $\{T_{\sigma}|\ \sigma\in
\mathfrak{S}_n\}$ forms a basis
for $H_n(0)$ and the multiplication of basis elements is
determined by:
$$T_iT_{\sigma}=\left\{\begin{array}{ll}
                  T_{s_i\sigma} & \hbox{ if } \ell(s_i\sigma)=\ell(\sigma)+1\\
                  -T_{\sigma}   & \hbox{ if }
                  \ell(s_i\sigma)=\ell(\sigma)-1.
                  \end{array}
                \right.$$
Here $\ell(\sigma)$ is the length of a reduced expression for
$\sigma$.

In $\mathfrak{S}_{m+n}$, we denote by $X_{(m,n)}$ the set of minimal
length coset representatives  of $\mathfrak{S}_{m+n}/\mathfrak{S}_m\times \mathfrak{S}_n$. We have
$\mathfrak{S}_{m+n}=\bigoplus_{\tau\in X_{(n,m)}} \tau(\mathfrak{S}_m\times \mathfrak{S}_n)$.
Moreover, our choice of representative implies that
$\ell(\tau\sigma)=\ell(\tau)+\ell(\sigma)$, for all $\tau\in X_{(n,m)}$ and   $\sigma\in
\mathfrak{S}_m\times \mathfrak{S}_n$~\cite{HJ}. This implies that
$$H_{m+n}(0)=\bigoplus_{\tau\in X_{(n,m)}}T_{\tau}(H_m(0)\otimes H_n(0)).$$ Therefore, when
we consider $H_{m+n}$ as a right $H_m(0)\otimes H_n(0)$-module it
is a direct sum of $(m+n)!/m!n!$ copies of $H_m(0)\otimes H_n(0)$.
Hence $H_{m+n}(0)$ is a right projective $H_m(0)\otimes
H_n(0)$-module. Analogously, $H_{m+n}(0)$ is a left projective
$H_m(0)\otimes H_n(0)$-module.

Now consider $H_N(0)$. To check the axiom (4) we need to better
understand the simple modules and projective indecomposable
modules of $H_N(0)$. For this we need to recall some results from
\cite{KT, PN}.
 For $i\in[1,N-1]$, let $\Box_i=1+T_i$.
These elements satisfy the relations
$$\begin{array}{rclc}
\Box_i^2&=&\Box_i,&\\
\Box_i\Box_j&=&\Box_j\Box_i& |i-j|>1,\\
\Box_i\Box_{i+1}\Box_i&=&\Box_{i+1}\Box_i\Box_{i+1}.&
\end{array}$$
In particular, the morphism defined by $T_i\longrightarrow\Box_i$
is an involution of $H_N(0)$. Since the $\Box_i$'s also satisfy
the braid relations, one can associate to each permutation
$\sigma\in \mathfrak{S}_N$ the element $\Box_{\sigma}$ of $H_N(0)$
defined by $\Box_{\sigma}=\Box_{i_1}\cdots\Box_{i_r}$, where
$s_{i_1}\cdots s_{i_r}$ is an arbitrary reduced decomposition of
$\sigma$.

For a composition $I=(i_1,\ldots,i_r)$ of $n$, the corresponding
\textit{ribbon diagram} of $I$ consists of $n$ boxes with $i_1$
boxes in the first row, $i_2$ boxes in the second row, $\cdots$,
$i_r$ boxes in the $r$th row and the first box in the next row is
under the last one in the previous row.
 We denote
by $\bar{I}=(i_r,\ldots,i_1)$ its mirror image and by $I\
\widetilde{}$\ \ its conjugate composition, ie., the composition
obtained by writing from right to left the lengths of the columns
of the ribbon diagram of $I$. For example, let $I=(3,1)$. Then
$\bar{I}=(1,3)$ and $I\ \widetilde{}=(2,1,1)$. The
corresponding ribbon diagrams are
$$\young(\ \ \ ,::\  )\hspace{1cm} \young(\  ,\ \ \ )
\hspace{1cm} \young( \ \ ,:\ ,:\ )$$

There are $2^{N-1}$ simple and $2^{N-1}$ indecomposable projective
$H_N$-modules. They can be realized as minimal left ideals and
indecomposable left ideals of $H_N(0)$ respectively. All the
simple modules are of dimension 1.

To describe the generators of the simple and indecomposable
projective $H_N$-modules, we associate with a composition $I$ of
$N$ two permutations $\alpha(I)$ and $\omega(I)$ of
$\mathfrak{S}_N$ defined by \begin{itemize}

\item $\alpha(I)$ is the permutation obtained by filling the
columns of the ribbon diagram of shape $I$ from bottom
to top and from left to right with the numbers $1,2,\ldots,N$;

\item $\omega(I)$ is the permutation obtained by filling the rows
of the ribbon diagram of shape $I$ from left to right
and from bottom to top with the numbers $1,2,\ldots,N$.
\end{itemize}

For example, consider the composition $I=(2,2,1,3)$ of $8$. The
fillings of the ribbon diagram of shape $I$ corresponding to
$\alpha(I)$ and $\omega(I)$ are
 $$ \begin{array}{ccc}
      \young(13,:26,::5,::478)& \ &
      \young(78,:56,::4,::123)\\
      \alpha(2,2,1,3)=13265478&\ &\omega(2,2,1,3)=78564123
 \end{array}$$

Let $I=(i_1,\ldots,i_r)$ be a composition and $\sigma\in
\mathfrak{S}_N$. The {\it descent set} of $\sigma$ is
$\hbox{Des}(\sigma)=\{i: \sigma(i)>\sigma(i+1)\}$ and we also
define $D(I)=\{i_1, i_1+i_2,\ldots,i_1+\cdots +i_{r-1}\}$. The
descent class $D_I=\{\sigma\in \mathfrak{S}_N:
\hbox{Des}(\sigma)=D(I)\}$ is the interval $[\alpha(I),\omega(I)]$
in the weak order on $\mathfrak{S}_N$ (see~\cite{KT} Lemma 5.2).

The simple $H_N(0)$-modules are indexed by all compositions of $N$.
The simple $H_N(0)$-module associated to a
composition $I$ is given by  the
minimal left ideal $C_I=H_N(0)\eta_I$, where
$\eta_I=T_{\omega(\bar{I})}\Box_{\alpha(I\ \widetilde{}\ )}$.
These modules form a complete system of simple $H_N(0)$-modules
and
$$T_i\eta_I=\left\{\begin{array}{cr}
                  -\eta_I & \hbox{ if }i\in D(I),\\
                  0   & \hbox{ if }i\notin D(I).
                  \end{array}
                \right.$$
We associates to $I$ an indecomposable projective
$H_N(0)$-module $M_I$ such that $M_I/\hbox{rad}(M_I)\cong H_N(0)\eta_I$.
This module is realized as the left ideal
$$M_I=H_N(0)\nu_I,$$
where $\nu_I=T_{\alpha(I)}\Box_{\alpha(\bar{I}\ \widetilde{}\ )}$.
A basis of $M_I$ is given by $\{T_{\sigma}\Box_{\alpha(\bar{I}\
\widetilde{}\ )}:\sigma\in [\alpha(I),\omega(I)]\}.$ The family
$(M_I)_{|I|=N}$ forms a complete system of projective
indecomposable $H_N(0)$-modules, and
$H_N(0)=\bigoplus_{|I|=N}M_I.$

\begin{remark} The results above are remarkable, specially considering the fact that
the $\nu_I$ are not the minimal idempotents of $H_N(0)$; this is an open problem in general.
The $\nu_I$ are not even idempotent in general. For example, let
$I=(2,1)$. Then $\bar{I}\ \widetilde{}=(1,2)$, $\alpha(I)=132=s_2$
and $\alpha(\bar{I}\ \widetilde{}\ )=213=s_1$.
$\nu_I^2=T_2\Box_1T_2\Box_1=T_2(1+T_1)T_2(1+T_1)=(T_2+T_2T_1)(T_2+T_2T_1)
=T_2^2+T_2^2T_1+T_2T_1T_2+T_2T_1T_2T_1=-T_2-T_2T_1+T_2T_1T_2-T_2T_1T_2
=-T_2-T_2T_1=-T_2(1+T_1)=-T_2\Box_1=-\nu_I\neq\nu_I.$
\end{remark}

{}From~\cite{PN}, we know that
$H_N(0)T_{\alpha(I)}\Box_{\alpha(\bar{I}\ \widetilde{}\ )}\cong
H_N(0)\Box_{\alpha(\bar{I}\ \widetilde{}\ )}T_{\alpha(I)}$ as left
ideals (also as left modules). Denote by  ``$-1$'' the anti-morphism
of $H_N(0)$ which reverses the order of the product of the
generators in all monomials.  For instance, $(T_{i_1}\cdots
T_{i_r})^{-1}=T_{i_r}\cdots T_{i_1}$, i.e.,
$(T_{\sigma})^{-1}=T_{\sigma^{-1}}$. This identity also holds when
we replace $T_i$ by $\Box_i$. Since $\alpha(I)^{-1}=\alpha(I)$ we
have $H_N(0)\nu_I\cong H_N(0)\nu_I^{-1}.$ Similarly,
$\nu_IH_N(0)\cong \nu_I^{-1}H_N(0)$ as right modules.

Let $g_I$ be a primitive idempotent such that $H_N(0)g_I\cong
H_N(0)\nu_I$. Obviously $g_I^{-1}$ is also a primitive idempotent
in $H_N(0)$ with $g_I^{-1}H_N(0)\cong\nu_I^{-1}H_N(0)\cong
\nu_IH_N(0)$. If $H_N(0)g_I\cong \bigoplus(H_k(0)\ot
H_l(0))(e_J\ot f_L)$ where $k+l=N$, $e_J$ and $f_L$ are primitive
idempotents in $H_k(0)$ and $H_l(0)$ respectively, then at the
same time we have $g_I^{-1}H_N(0)\cong \bigoplus(e_J^{-1}\ot
f_L^{-1})(H_k(0)\ot H_l(0))$. To show axiom (4) we need an
auxiliary result:

\begin{proposition}  Let $H$ be a self-injective
algebra and $g$ be an element in $H$ such that $Hg$ is a
projective $H$-module. Then $Hg\cong H\nu$ as $H$-modules for some
$\nu\in H$ if and only if there exist $a,b,c,d\in H$ such that
$a\nu=gb,\ cg=\nu b,\ acg=g,\ ca\nu=\nu,\ gbd=g$ and $\nu db=\nu.$
\end{proposition}

\begin{proof} Suppose that there exist $a,b,c,d\in
H$ such that $a\nu=gb,\ cg=\nu b,\ acg=g,\ ca\nu=\nu,\ gbd=g$ and
$\nu db=\nu.$ Define $\phi:Hg\ra H\nu$ as a (left) $H$-module
homomorphism by $\phi(g)=a\nu$. Then
$\phi(cg)=c\phi(g)=ca\nu=\nu$. Define $\psi:H\nu\ra Hg$ as a
(left) $H$-module homomorphism by $\psi(\nu)=cg$. Since
$(\phi\circ\psi)(\nu)=\phi(\psi(\nu))=\phi(cg)=c\phi(g)=ca\nu=\nu$
and
$(\psi\circ\phi)(g)=\psi(\phi(g))=\psi(a\nu)=a\psi(\nu)=acg=g$,
$\psi=\phi^{-1}$ and $\phi$ is an isomorphism from $Hg$ to $H\nu$.

Conversely, suppose that $H$ is a self-injective algebra, $g$ is
an element in $H$ such that $Hg$ is a projective $H$-module. Let
$\phi:Hg\ra H\nu$ be a (left) $H$-module isomorphism. Then
$\phi(g)=a\nu$ and $\phi^{-1}(\nu)=cg$ for some $a,\ c\in H.$
Hence $\nu=\phi(cg)=c\phi(g)=ca\nu$ and
$g=\phi^{-1}(a\nu)=a\phi^{-1}(\nu)=acg.$ Since $H$ is
self-injective, i.e., an  $H$-module is projective if and only if
it is injective~\cite{ARS}, $H\nu$ is an injective module and
$\phi:Hg\ra H\nu$ can be extended to a homomorphism from $H$ to
$H\nu$ such that the following diagram
    \begin{diagram}
     Hg & \rInto& H\\
    \dTo^{\phi} &\ldDashto_{\exists!\phi} \\
    H\nu
\end{diagram}
is commutative. For convenience we also write the homomorphism
$\phi:H\ra H\nu$. Similarly, $\phi^{-1}:H\nu\ra Hg$ can be
extended to a homomorphism $\phi^{-1}:H\ra Hg$. Let $\phi(1)=b$
and $\phi^{-1}(1)=d$ for some $b,\ d\in H$. Then
$a\nu=\phi(g)=g\phi(1)=gb$, $cg=\phi^{-1}(\nu)=\nu\phi^{-1}(1)=\nu
d$, $\nu=\phi(cg)=\phi(\nu d)=\nu d\phi(1)=\nu db$ and
$g=\phi^{-1}(a\nu)=\phi^{-1}(gb)=gb\phi^{-1}(1)=gbd.$
\end{proof}

Since $H_N(0)$ is self-injective~\cite{MF}, we have

\begin{corollary} If $H_N(0)g_I\cong H_N(0)\nu_I$ for some
primitive idempotent $g_I\in H_N(0)$, then $H_N(0)g_I^{-1}\cong
H_N(0)\nu_I$, i.e., $H_N(0)g_I\cong H_N(0)g_I^{-1}.$  Similarly
$g_IH_N(0)\cong g_I^{-1}H_N(0).$
\end{corollary}

\begin{proof}
Since $H_N(0)g_I\cong H_N(0)\nu_I$ there exist $a,b,c,d\in H_N(0)$
such that $a\nu_I = g_Ib$, $cg_I = \nu_Id$, $acg_I = g_I$, $ca\nu_I = \nu_I,$, $g_Ibd = g_I$ and
$\nu_Idb =\nu_I$.
Applying ``${}^{-1}$'' to these equations we get
$d^{-1}\nu_I^{-1} = g_I^{-1}c^{-1}$, $ b^{-1}g_I^{-1} = \nu_I^{-1}a^{-1} $,
$  d^{-1}b^{-1}g_I^{-1} = g_I^{-1} $, $  b^{-1}d^{-1}\nu_I^{-1} = \nu_I^{-1} $,
$  g_I^{-1}c^{-1}a^{-1} = g_I^{-1} $ and $\nu_I^{-1}a^{-1}c^{-1} =\nu_I^{-1}$.
Setting $a'=d^{-1},b'=c^{-1},c'=b^{-1}$ and $d'=a^{-1}$ we obtain the equations
needed to show that $H_N(0)g_I^{-1}\cong H_N(0)\nu_I^{-1}\cong H_N(0)\nu_I$.
\end{proof}

Hence, condition (4) holds.

Next we prove the identity in condition (5) for $G_0(A)$. First we
need to introduce the definition of shuffle. Let $A$ be a totally
ordered alphabet. $A^*$ denotes the set of all finite-length words
formed from the elements in $A$. The \textit{shuffle} is the
bilinear operation of $\mathbb{N}\langle A\rangle$~\cite{KT}
denoted by $\shuf$ and recursively defined on words by the
relations
$$\begin{array}{rcl}
1\shuf u&=&u\shuf 1=u,\\
(au)\shuf(bv)&=&a(u\shuf bv)+b(au\shuf v),
\end{array}$$where $1$ is the empty word, $u,v\in A^*$ and $a,b\in
A.$ One can show that $\shuf$ is associative. For convenience, we
also denote $u\shuf v$ the set of all words occur in the sum of
the shuffle. For instance,
$$21\shuf34=2134+2314+2341+3214+3241+3421.$$
It also means that $21\shuf34=\{2134,2314,2341,3214,3241,3421\}.$

{}From Proposition 5.7 in~\cite{KT}, let $I$ and $J$ be compositions
of $m$ and $n$. Let also $\sigma\in \mathfrak{S}_{[1,m]}$ and
$\tau\in\mathfrak{S}_{[m+1,m+n]}$ such that
$\hbox{Des}(\sigma)=D(I)$ and $\hbox{Des}(\tau)=D(J)$. Then
$$[\hbox{Ind}_{H_m(0)\ot H_n(0)}^{H_{m+n}(0)}C_I\ot C_J]
=\sum_{\omega \in \sigma\shuf\tau}[C_{C(\omega)}],$$ where
$C(\omega)$ denotes the composition associated with the descent
set of $\omega$.

\begin{proposition}The following identity holds
     $$\begin{array}{l}
          a_1\cdots a_m\shuf b_1\cdots b_n\\
          \qquad\quad=
          \sum_{i=0}^k(a_1\cdots a_i\shuf b_1\cdots b_{k-i})(a_{i+1}\cdots a_m\shuf b_{k-i+1}\cdots b_n)
         \end{array}$$
for any $0\leq k\leq m$.
\end{proposition}

\begin{proof} We proceed by induction.
When $k=0$, we have the trivial identity
$$a_1\cdots a_m\shuf b_1\cdots b_n=(1\shuf1)(a_1\cdots a_m\shuf b_1\cdots b_n).$$
For  $k=1$, we obtain the defining recursion of shuffle:
     $$\begin{array}{l}
     a_1\cdots a_m\shuf b_1\cdots b_n\\
     \qquad =
         (1\shuf b_1)(a_1\cdots a_m\shuf b_2\cdots b_n)
         +(a\shuf1)(a_2\cdots a_m\shuf b_1\cdots b_n).
      \end{array}$$

For $k>1$ we start with
  $$a_1\cdots a_m\shuf b_1\cdots b_n = a_1(a_2\cdots a_m\shuf b_1\cdots b_n)+
     b_1(a_1\cdots a_m\shuf b_2\cdots b_n)$$
and use the induction hypothesis to get
  $$\begin{array}{l}
   a_1(a_2\cdots a_m\shuf b_1\cdots b_n)+b_1(a_1\cdots a_m\shuf b_2\cdots b_n)\\
   \qquad = a_1\sum_{i=1}^{k+1} (a_2\cdots a_i\shuf b_1\cdots
             b_{k+1-i})(a_{i+1}\cdots a_i\shuf b_{k+1-i+1}\cdots b_n)\\
             \qquad \quad +b_1\sum_{i=0}^k (a_1\cdots a_i\shuf b_2\cdots
             b_{k+1-i})(a_{i+1}\cdots a_i\shuf b_{k+1-i+1}\cdots b_n)\\
   \qquad = b_1(1\shuf b_2\cdots b_{k+1})(a_1\cdots a_m\shuf b_{k+2}\cdots b_n)\\
             \qquad \quad+\sum_{i=1}^k a_1(a_2\cdots a_i\shuf b_1\cdots
        b_{k+1-i})(a_{i+1}\cdots a_i\shuf b_{k+1-i+1}\cdots b_n)\\
             \qquad \quad+\sum_{i=1}^k b_1(a_1\cdots a_i\shuf b_2\cdots
        b_{k+1-i})(a_{i+1}\cdots a_i\shuf b_{k+1-i+1}\cdots b_n)\\
             \qquad \quad+a_1(a_2\cdots a_{k+1}\shuf 1)(a_{k+2}\cdots a_m\shuf b_{1}\cdots b_n)\\
   \qquad =(1\shuf b_1\cdots b_{k+1})(a_1\cdots a_m\shuf b_{k+2}\cdots b_n)\\
             \qquad \quad+\sum_{i=1}^k(a_1\cdots a_i\shuf b_1\cdots
        b_{k+1-i})(a_{i+1}\cdots a_i\shuf b_{k+1-i+1}\cdots b_n)\\
             \qquad \quad+(a_1\cdots a_{k+1}\shuf1)(a_{k+2}\cdots a_m\shuf b_{1}\cdots b_n)\\
   \qquad =\sum_{i=0}^{k+1}(a_1\cdots a_i\shuf b_1\cdots
        b_{k+1-i})(a_{i+1}\cdots a_i\shuf b_{k+1-i+1}\cdots b_n).\\
 \end{array}$$
\end{proof}

This implies that condition (5) holds for $G_0(A)$.
\end{example}

\begin{remark}
Consider the direct sum of 0-Hecke-Clifford algebras~\cite{BHT}
$HCl_n(0),\ n\geq0,$ which are superalgebras. Here again it is possible to check that this tower satisfies the
modified conditions (1)-(5) to show that  $\bigoplus_{n\geq0}HCl_n(0)$ is
also a tower of superalgebras.
\end{remark}

\subsection{An example not satisfying Condition (5)}

If one considers a direct sum of algebras that does not satisfy
condition (3) then the induction and restriction may not be well
defined. If it does not satisfy condition (4), then its
Grothendieck groups are graded Hopf algebras respectively but not
necessarily dual to each other. Hence we are mostly interested in finding structure that satisfies
all our axioms but (5). We give some in~\cite{BL} but the simplest one was given to us by F.~Hivert:

\begin{example}
Let $A_n=\mathbb{C}[\mathbb{Z}/2\mathbb{Z}]^{\ot n}$ and
$\rho_{m,n}:\mathbb{C}[\mathbb{Z}/2\mathbb{Z}]^{\ot m}\ot
\mathbb{C}[\mathbb{Z}/2\mathbb{Z}]^{\ot n}\ra
\mathbb{C}[\mathbb{Z}/2\mathbb{Z}]^{\ot (m+n)}$ be the identity
map. It is clear that this tower satisfies all conditions (1)-(4).
It does not satisfy Condition (5). To see this, we know that there
are two simple $A_1$-modules $T$, the trivial module and $S$, the
sign module. They are also indecomposable projective
$A_1$-modules. Any simple (or indecomposable projective)
$A_n$-module is an $n$-tensor product of $T$'s and $S$'s. To see
that (5) is not satisfied in general, consider the left hand side
of the formula
$$[\hbox{Res}_{A_1\ot A_1}^{A_2}\hbox{Ind}_{A_1\ot A_1}^{A_2}(T\ot S)]
 =[\hbox{Res}_{A_1\ot A_1}^{A_2}(T\ot S)]
=[T\ot S].$$
But the right hand side is
$$\begin{array}{l}
[\widetilde{\hbox{Ind}}_{A_0\ot A_1\ot A_1\ot A_0}^{A_1\ot A_1}
(\hbox{Res}_{A_0\ot A_1}^{A_1}T\ot \hbox{Res}_{A_1\ot
A_0}^{A_1}S)]\\\quad  + \quad [\widetilde{\hbox{Ind}}_{A_1\ot A_0\ot A_0\ot A_1}^{A_1\ot
A_1}
(\hbox{Res}_{A_1\ot A_0}^{A_1}T\ot \hbox{Res}_{A_0\ot A_1}^{A_1}S)]  = [S\ot T]+[T\ot S].\\
\end{array}$$
These are not equal.
\end{example}

\begin{remark}
The algebra $A=(\bigoplus A_n,\rho)$ above does not satisfy our
condition (5) and its Grothendieck groups $G(A)$ and $K(A)$ are
not Hopf algebra in the strict sense. Yet, in this case  $G_0(A)$
and $K_0(A)$ are generalized bialgebras in the sense of
Loday~\cite{L}. the multiplication $\pi$ and the comultiplication
$\Delta$ satisfies a very simple compatibility relation. Let
$\hat\Delta(x)=\Delta(x)-1\ot x-x\ot 1$. Then
  \begin{eqnarray}\label{AsAs}
  \hat\Delta\circ\pi = Id \ot Id + (\pi\ot Id)\circ(Id\ot\hat\Delta) + (Id\ot\pi)\circ(\hat\Delta\ot\pi).
  \end{eqnarray}
At the module level, this is equivalent to the following requirement:

 \noindent
(5)'  In $G_0(A)$ we have
 $$\begin{array}{l}
         [\mbox{Res}_{A_k\ot A_{m+n-k}}^{A_{m+n}}\mbox{Ind}_{A_m\ot A_n}^{A_{m+n}}(M\ot N)]\\
           \qquad\quad =
             \left\{\begin{array}{lr}
               \big[(Id\ot\mbox{Ind}_{A_{m-k}\ot A_{n}}^{A_{m+n-k}})
                    ((\mbox{Res}_{A_k\ot A_{m-k}}^{A_{m}}M)\ot N)\big]
                         & \hbox{ if }k<m,\\
                 \big[M\ot N\big]   & \hbox{ if }k=m,\\
                \big[(\mbox{Ind}_{A_{m}\ot A_{k-m}}^{A_{k}}\ot Id)
                    (M\ot (\mbox{Res}_{A_{k-m}\ot A_{m+n-k}}^{A_{n}}M))\big]
                         & \hbox{ if }k>m.\\
                  \end{array}
                \right.
           \end{array}$$
 This is easy to check for  $G_0(A)$. It is thus a self-dual bialgebra
 satisfying the compatibility relation~(\ref{AsAs}).
\end{remark}

\section{Concluding remarks}

In the last example of Section~4, we encountered a graded algebra that
satisfies our conditions (1)--(4) but not (5). Yet, following Loday~\cite{L},
we still have an interesting (generalized) bialgebra structure  on its Grothendieck groups.
We have given an alternative axiom, (5)', that shows that we get the kind of algebra
satisfying the compatibility relations ~(\ref{AsAs}).

This open the door to many avenues. The conditions (1)--(4) on a
graded algebra $A$ are essential to make sure that we can define a
structure of graded algebra and of graded coalgebra on $G_0(A)$
and $K_0(A)$ with duality. Then one may ask what kind of
compatibility one can get between the induction and the
restriction. In this sense there are many alternatives to our
condition (5). It would be interesting to find what is the
required condition for each of the generalized bialgebras
of~\cite{L} and to give examples for each cases. One can also
define different kinds of inductions and restrictions to allow for
different kind of operations on the Grothendieck groups of the
tower. This is left to future work.


\end{document}